\documentclass[11pt]{article}
\usepackage{epsfig}

\def\be{\begin{equation}}
\def\ee{\end{equation}}
\def\bea{\begin{eqnarray}}
\def\eea{\end{eqnarray}}
\def\bes{\begin{eqnarray*}}
\def\ees{\end{eqnarray*}}

\def\nn{\nonumber}
\def\lb{\label}
\def\bs{\setminus}
\def\pt{\partial}

\def\R{{\bf R}}
\def\C{{\bf C}}
\def\Z{{\bf Z}}

\def\N{{\bf N}}

\def\Q{{\bf Q}}

\def\RP{{\bf RP}}

\def\ga{{\gamma}}

\def\th{{\theta}}

\def\Om{{\Omega}}

\def\lm{{\lambda}}
\def\Lm{{\Lambda}}

\def\rank{{\rm rank}}
\def\Sp{{\rm Sp}}

\def\ol{\overline}

\def\hb{\vrule height0.18cm width0.14cm $\,$}

\title{Two closed geodesics on  compact  bumpy Finsler manifolds}

\author{
Wei Wang\thanks{Partially supported by NSFC No. 11222105, 11431001, E-mail: wangwei@math.pku.edu.cn}\\\\
 School of Mathematical Sciences and LMAM \\Peking University, Beijing 100871\\
The People's Republic of China\\}

\begin{document}

\maketitle

\begin{abstract}
{\it In this paper, we prove there are at least two closed geodesics
 on any compact  bumpy Finsler $n$-manifold with finite fundamental group and $n\ge 2$.
 Thus generically there  are at least two closed geodesics
 on  compact  Finsler manifolds with finite fundamental group. 
Furthermore,  there are at least two closed geodesics
 on any compact  Finsler $2$-manifold,
and this lower bound is achieved
 by the Katok 2-sphere $(S^2, F)$ and 2-real projective space $(S^2/\Z_2, F)$, cf.  \cite{Kat}. 

 }
\end{abstract}

{\bf Key words}: Closed geodesic, Finsler manifold, bumpy.

{\bf 2000 Mathematics Subject Classification}: 53C22, 58E05, 58E10.

\renewcommand{\theequation}{\thesection.\arabic{equation}}
\renewcommand{\thefigure}{\thesection.\arabic{figure}}

\setcounter{equation}{0}
\section{Introduction and main results}

It is well-known that there are at least two closed geodesics on all compact
bumpy Riemannian manifolds $M$ with ${\rm dim} M\ge2$ except for some special 
type of manifolds, cf. \cite{Fet} or Theorem 4.1.8 of \cite{Kli2} and Remark 1.5 below. 
While its proof
depends on the symmetric property for the Riemannian metric, and consequently the proof
 carries over to the symmetric Finsler case. But for non-symmetric Finsler case, the proof does not work, hence we must develop new methods to handle the problem in this case. This paper is devoted to do this.

 Let us recall firstly the definition of the Finsler
metrics.

{\bf Definition 1.1.} (cf. \cite{BCS} or \cite{She}) {\it Let $M$ be a finite
dimensional smooth manifold. A function $F:TM\to [0,+\infty)$ is a {\rm
Finsler metric} if it satisfies

(F1) $F$ is $C^{\infty}$ on $TM\bs\{0\}$,

(F2) $F(x,\lm y) = \lm F(x,y)$ for all $x\in M$, $y\in T_xM$ and
$\lm>0$,

(F3) For every $y\in T_xM\bs\{0\}$, the quadratic form
$$ g_{x,y}(u,v) \equiv
         \frac{1}{2}\frac{\pt^2}{\pt s\pt t}F^2(x,y+su+tv)|_{t=s=0},
         \qquad \forall u, v\in T_xM, $$
is positive definite.

In this case, $(M,F)$ is called a {\rm Finsler manifold}. $F$ is
{\rm symmetric} if $F(x,-y)=F(x,y)$ holds for all $x\in M$ and
$y\in T_xM$. $F$ is {\rm Riemannian} if $F(x,y)^2=\frac{1}{2}G(x)y\cdot
y$ for some symmetric positive definite matrix function $G(x)\in
GL(T_xM)$ depending on $x\in M$ smoothly. }

A closed curve on a Finsler manifold is a {\it closed geodesic} if it is
locally the shortest path connecting any two nearby points on this
curve (cf. \cite{She}). As usual, on any Finsler manifold
$(M, F)$, a closed geodesic $c:S^1=\R/\Z\to M$ is {\it prime}
if it is not a multiple covering (i.e., iteration) of any other
closed geodesics. Here the $m$-th iteration $c^m$ of $c$ is defined
by $c^m(t)=c(mt)$, where $m\in\N$. The inverse curve $c^{-1}$ of $c$ is defined by
$c^{-1}(t)=c(1-t)$ for $t\in \R$.  Note that unlike Riemannian manifold,
the inverse curve $c^{-1}$ of a closed geodesic $c$
on a non-symmetric Finsler manifold need not be a geodesic.
Two prime closed geodesics
$c$ and $d$ are {\it distinct} if there is no $\th\in (0,1)$ such that
$c(t)=d(t+\th)$ for all $t\in\R$.
We shall omit the word {\it distinct} when we talk about more than one prime closed geodesic.
On a symmetric Finsler (or Riemannian) manifold, two closed geodesics
$c$ and $d$ are called { \it geometrically distinct} if $
c(S^1)\neq d(S^1)$, i.e., their image sets in $M$ are distinct.

For a closed geodesic $c$ on $(M,\,F)$, denote by $P_c$
the linearized Poincar\'{e} map of $c$. Then $P_c\in \Sp(2n-2)$ is symplectic.
A closed geodesic $c$ is called
{\it non-degenerate} if $1$
is not an eigenvalue of $P_c$. A Finsler manifold $(M,\,F)$
is called {\it bumpy} if all the closed geodesics on it are
non-degenerate. Note that  bumpy Finsler metrics are generic in the
set of Finsler metrics, cf. \cite{Rad4}.

The following are the main results in this paper:

{\bf Theorem 1.2.} {\it There exist at least two prime closed geodesics on every compact  bumpy Finsler manifold $(M, F)$ with finite fundamental group and $\dim M\ge 2$. }

Furthermore, if $\dim M=2$, the bumpy and finite fundamental group conditions are not needed, and we have the following:

{\bf Theorem 1.3.} {\it There exist at least two prime closed geodesics on every compact   Finsler manifold $(M, F)$ with $\dim M=2$. }

{\bf Remark 1.4. } In 1973,  Katok  in \cite{Kat}  found some non-symmetric 
Finsler metrics on CROSSs (compact rank one symmetric spaces) with only finitely many prime closed geodesics and all closed
geodesics are non-degenerate. The  number of closed geodesics
 on $S^n$ that one obtains in these examples is $2[\frac{n+1}{2}]$, 
where $[a]=\max\{k\in\Z\,|\,k\le a\}$ for $a\in\R$, cf. \cite{Zil}.
 
 We are aware of a number of results concerning closed geodesics on Finsler manifolds. According to the classical theorem of Lyusternik-Fet
\cite{LyF} from 1951, there exists at least one closed geodesic on
every compact Riemannian manifold. The proof of this theorem is
variational and carries over to the Finsler case.  In \cite{BaL}, V. Bangert and Y. Long proved that on any
Finsler 2-sphere $(S^2, F)$, there exist at least two  closed
geodesics. In \cite{Rad3}, H.-B. Rademacher studied the existence and
stability of closed geodesics on positively curved Finsler
manifolds.
In \cite{DuL1} of Duan and Long and in  \cite{Rad4} of Rademacher, they proved there exist 
at least two  closed geodesics on any bumpy
Finsler $n$-sphere independently. In  \cite{Rad5}, Rademacher proved there exist at least two  closed geodesics on any bumpy Finsler $\bf{CP^2}$.
In \cite{DLW}, Duan, Long and Wang proved there exist at least two 
closed geodesics on any compact simply-connected bumpy Finsler manifold. 
In \cite{DLX}, Duan, Long and Xiao  proved the existence
of at least two  non-contractible closed geodesics on any bumpy Finsler $\RP^3$. 
In \cite{Tai2}, Taimanov  proved the existence of at least two  non-contractible closed geodesics
on any bumpy Finsler $\RP^2$. In \cite{LiX}, Liu and Xiao proved there exist at least two non-contractible  closed geodesics on
any bumpy Finsler ${\rm RP}^n$.

{\bf Remark 1.5.}  For the case $\pi_1(M)$ is infinite,
as pointed out in \cite{Tai1} of  I. Taimanov, there are  at least two prime closed geodesics on all compact   Riemannian  manifolds  except  $M$ is  an Eilenberg-MacLane complex 
$K(\pi, 1)$  such that  $\pi$ is different from $\Z$ and contains an element $g\in\pi$
such that any element of $\pi$ is conjugate to one power of $g$.
cf. Theorem 4.1.8 of \cite{Kli2} and \cite{BaH} for a proof.  Note that the same  proof yields 
there are  at least two prime closed geodesics on all compact   Finsler  manifolds  with infinite
fundamental group except  $M$ is  an Eilenberg-MacLane complex $K(\pi, 1)$ as above.

{\bf Remark 1.6.} In \cite{LLX}, Liu, Long and Xiao
proved  in every non-trivial homotopy class $\alpha\in\pi_1(S^n/\Gamma)$  
with finite order of a bumpy Finsler $S^n/\Gamma$,
where $\Gamma$ is a finite group acts on $S^n$ freely  and isometrically,
there exist at least two distinct closed geodesics.  While 
we don't know whether the two closed geodesics obtained in Theorem 1.2
belongs to some identical homotopy class. Note also that the result of \cite{LLX} does not imply
there are more than two  prime closed geodesics on $S^n/\Gamma$. For example, 
let $(S^2, F)$ be the Katok 2-sphere, cf. \cite{Kat} or \cite{Zil}, then $(S^2/\Z_2, F)$
has exactly two closed geodesics $c_+$ and $c_-$ in $0\neq \alpha\in\pi_1(S^2/\Z_2)$,
and exactly two closed geodesics $c_+^2$ and $c_-^2$ in $0=\alpha^2\in\pi_1(S^2/\Z_2)$.

In this paper, let $\N$, $\N_0$, $\Z$, $\Q$, $\R$, and $\C$ denote
the sets of natural integers, non-negative integers, integers,
rational numbers, real numbers, and complex numbers respectively.
We use only singular homology modules with $\Q$-coefficients.
For terminologies in algebraic topology we refer to \cite{GrH}.
For $k\in\N$, we denote by $\Q^k$ the direct sum $\Q\oplus\cdots\oplus\Q$ of
$k$ copies of $\Q$ and $\Q^0=0$. For an $S^1$-space $X$,
we denote by $\overline{X}$ the quotient space $X/S^1$.

\setcounter{equation}{0}
\section{Critical point theory for closed geodesics}

Let $M=(M,F)$ be a compact Finsler manifold, the space
$\Lambda=\Lambda M$ of $H^1$-maps $\gamma:S^1\rightarrow M$ has a
natural structure of Riemannian Hilbert manifolds on which the group
$S^1=\R/\Z$ acts continuously by isometries, c.f. \cite{Kli1}-\cite{Kli3}. This action is defined
by $(s\cdot\gamma)(t)=\gamma(t+s)$ for all $\gamma\in\Lm$ and $s,
t\in S^1$. For any $\gamma\in\Lambda$, the energy functional is
defined by
\be E(\gamma)=\frac{1}{2}\int_{S^1}F(\gamma(t),\dot{\gamma}(t))^2dt.
\lb{2.1}\ee
It is $C^{1,1}$ and invariant under the $S^1$-action,  cf. \cite{Mer}.
The critical points of $E$ of positive energies are precisely the
closed geodesics $\gamma:S^1\to M$. The index form of the functional
$E$ is well defined along any closed geodesic $c$ on $M$, which we
denote by $E''(c)$. As usual we define the index $i(c)$ of $c$ as the maximal dimension of
subspaces of $T_c \Lambda$ on which $E^{\prime\prime}(c)$ is negative definite, and the
nullity $\nu(c)$ of $c$ so that $\nu(c)+1$ is the dimension of the null
space of $E^{\prime\prime}(c)$, cf. Definition 2.5.4 of \cite{Kli3}.
 In the following, we denote
by
\be \Lm^\kappa=\{d\in \Lm\;|\;E(d)\le\kappa\},\quad
\Lm^{\kappa-}=\{d\in \Lm\;|\; E(d)<\kappa\},
  \quad \forall \kappa\ge 0. \lb{2.2}\ee
For a closed geodesic $c$ we set $ \Lm(c)=\{\ga\in\Lm\mid E(\ga)<E(c)\}$.

For $m\in\N$ we denote the $m$-fold iteration map
$\phi_m:\Lambda\rightarrow\Lambda$ by $\phi_m(\ga)(t)=\ga(mt)$, for all
$\,\ga\in\Lm, t\in S^1$, as well as $\ga^m=\phi_m(\gamma)$. If $\gamma\in\Lambda$
is not constant then the multiplicity $m(\gamma)$ of $\gamma$ is the order of the
isotropy group $\{s\in S^1\mid s\cdot\gamma=\gamma\}$. For a closed geodesic $c$,
the mean index $\hat{i}(c)$ is defined as usual by
$\hat{i}(c)=\lim_{m\to\infty}i(c^m)/m$. 

We call a closed geodesic satisfying the isolation condition, if
the following holds:

{\bf (Iso)  For all $m\in\N$ the orbit $S^1\cdot c^m$ is an
isolated critical orbit of $E$. }

Note that if the number of prime closed geodesics on a Finsler manifold
is finite, then all the closed geodesics satisfy (Iso).

Using singular homology with rational
coefficients we consider the following critical $\Q$-module of a closed geodesic
$c\in\Lambda$:
\be \overline{C}_*(E,c)
   = H_*\left((\Lm(c)\cup S^1\cdot c)/S^1,\Lm(c)/S^1\right). \lb{2.3}\ee

{\bf Proposition 2.1.} (cf. Satz 6.11 of \cite{Rad2} or Proposition
3.12 of \cite{BaL}) {\it Let $c$ be a prime closed geodesic on a Finsler
manifold $(M,F)$ satisfying (Iso). Then we have
\bea \overline{C}_q( E,c^m)
&&\equiv H_q\left((\Lm(c^m)\cup S^1\cdot c^m)/S^1, \Lm(c^m)/S^1\right)\nn\\
&&= \left(H_{i(c^m)}(U_{c^m}^-\cup\{c^m\},U_{c^m}^-)
    \otimes H_{q-i(c^m)}(N_{c^m}^-\cup\{c^m\},N_{c^m}^-)\right)^{+\Z_m} \nn
\eea

(i) When $\nu(c^m)=0$, there holds
\bea \overline{C}_q( E,c^m) = \left\{\matrix{
     \Q, &\quad {\it if}\;\; i(c^m)-i(c)\in 2\Z,\;{\it and}\;
                   q=i(c^m),\;  \cr
     0, &\quad {\it otherwise}. \cr}\right.  \nn\eea

(ii) When $\nu(c^m)>0$, there holds
\bea \overline{C}_q( E,c^m) =
H_{q-i(c^m)}(N_{c^m}^-\cup\{c^m\}, N_{c^m}^-)^{(-1)^{i(c^m)-i(c)}\Z_m},\nn\eea
where $N_{c^m}$ is a local characteristic manifold at $c^m$ and $N^-_{c^m}=N_{c^m}\cap \Lambda(c^m)$,  $U_{c^m}$ 
is a local negative disk at $c^m$ and $U^-_{c^m}=U_{c^m}\cap \Lambda(c^m)$,
$H_{\ast}(X,A)^{\pm\Z_m}
   = \{[\xi]\in H_{\ast}(X,A)\,|\,T_{\ast}[\xi]=\pm [\xi]\}$
   where $T$ is a generator of the $\Z_m$-action.}

{\bf Definition 2.2.} {\it The Euler characteristic $\chi(c^m)$
of $c^m$ is defined by
\bea \chi(c^m)
&\equiv& \chi\left((\Lm(c^m)\cup S^1\cdot c^m)/S^1, \Lm(c^m)/S^1\right), \nn\\
&\equiv& \sum_{q=0}^{\infty}(-1)^q\dim \overline{C}_q( E,c^m).
\lb{2.4}\eea
Here $\chi(A, B)$ denotes the usual Euler characteristic of the space pair $(A, B)$.

The average Euler characteristic $\hat\chi(c)$ of $c$ is defined by }
\be \hat{\chi}(c)=\lim_{N\to\infty}\frac{1}{N}\sum_{1\le m\le N}\chi(c^m).
\lb{2.5}\ee
By Remark 5.4 of \cite{Wan},  $\hat\chi(c)$ is well-defined and is a
rational number. In particular, if $c^m$ are non-degenerate for $\forall m\in\N$, then 
 \bea
 \hat\chi(c)=\;\;\left\{\matrix{
    (-1)^{i(c)},&\quad {\it if}\quad i(c^2)-i(c)\in2\Z,  \cr
    \frac{(-1)^{i(c)}}{2},&\quad {\it otherwise.}\cr
    }\right.  \lb{2.6}\eea

Set $\ol{\Lm}^0=\ol{\Lambda}^0M =\{{\rm
constant\;point\;curves\;in\;}M\}\cong M$. Let $(X,Y)$ be a
space pair such that the Betti numbers $b_i=b_i(X,Y)=\dim
H_i(X,Y;\Q)$ are finite for all $i\in \Z$. As usual the {\it
Poincar\'e series} of $(X,Y)$ is defined by the formal power series
$P(X, Y)=\sum_{i=0}^{\infty}b_it^i$. We need the following  results on Betti numbers.

For a compact and simply-connected Finsler manifold $M$ with
$H^*(M;\Q)\cong T_{d,h+1}(x)$ with the generator $x$ of degree $d$ and height $h + 1$,
if  $d$ is odd, then $x^2=0$ and
$h=1$ in $T_{d,h+1}(x)$, thus $M$ is rationally homotopy equivalent to
$S^d$ (cf.  \cite{Rad1} or \cite{Hin}).

\medskip

{\bf Proposition 2.3.} (cf.  Theorem 2.4, Remark 2.5 of \cite{Rad1} and Lemma 2.5, 2.6 of \cite{DuL2})
{\it Let $M$ be a compact
simply-connected manifold with $H^*(M;\Q)\cong T_{d,h+1}(x)$.
Then the Betti numbers of the free loop space of $M$ defined by
$b_q = \rank H_q(\Lm M/S^1,\Lm^0 M/S^1;\Q)$ for $q\in\Z$ are given by

(i) If $h=1$ and $d\in 2\N+1$,  then we have
\bea {b}_q = \;\;\left\{\matrix{
    2,&\quad {\it if}\quad q\in {\cal K}\equiv\{k(d-1)\,|\,2\le k\in\N\},  \cr
    1,&\quad {\it if}\quad q\in \{d-1+2k\,|\,k\in\N_0\}\setminus {\cal K},\cr
        0 &\quad {\it otherwise}. \cr}\right. \lb{2.7}\eea

(ii) If $h=1$ and $d\in 2\N$,  then we have
\bea {b}_q = \;\;\left\{\matrix{
    2,&\quad {\it if}\quad q\in {\cal K}\equiv\{k(d-1)\,|\,3\le k\in(2\N+1)\},  \cr
    1,&\quad {\it if}\quad q\in \{d-1+2k\,|\,k\in\N_0\}\setminus {\cal K},\cr
        0 &\quad {\it otherwise}. \cr}\right. \lb{2.8}\eea

(iii) If $h\ge 2$ and $d\in 2\N$. Let $D=d(h+1)-2$ and
\bea \Om(d,h) = \{k\in 2\N-1&\,|\,& iD\le k-(d-1)=iD+jd\le iD+(h-1)d\;  \nn\\
         && \mbox{for some}\;i\in\N\;\mbox{and}\;j\in [1,h-1]\}. \lb{2.9}\eea
Then we have
\bea b_q = \left\{\matrix{
    0, & \quad \mbox{if}\ q\in2\Z \mbox{ or}\ q\le d-2,  \cr
    [\frac{q-(d-1)}{d}]+1, & \quad \mbox{if}\ q\in 2\N-1\;\mbox{and}\;d-1\le q < d-1+(h-1)d, \cr
    h+1, & \quad \mbox{if}\ q\in \Om(d,h), \cr
    h, & \quad \mbox{otherwise}. \cr}\right.
\lb{2.10}\eea}

By a similar proof of Theorem 5.5  of \cite{Wan}, we have the following mean index identity:

{\bf Proposition 2.4.} (cf. Theorem 3.1 of \cite{Rad1} and Satz 7.9 of \cite{Rad2}) {\it Let
$(M,F)$ be a compact simply-connected  Finsler manifold with
$\,H^{\ast}(M,\Q)=T_{d,h+1}(x)$ and possess finitely many prime closed
geodesics. Denote the prime closed geodesics on $(M,F)$
with positive mean indices by $\{c_j\}_{1\le j\le q}$ for some $q\in\N$.
Then the following identity holds
\be \sum_{j=1}^q\frac{\hat\chi(c_j)}{\hat{i}(c_j)}=B(d,h)
=\left\{\matrix{
     -\frac{h(h+1)d}{2d(h+1)-4}, &\quad d\;\;{\it is\;even},\cr
     \frac{d+1}{2d-2}, &\quad d\;\;{\it is\;odd}\;(then \; h=1),\cr}\right.   \lb{2.11}\ee
where $\dim M=hd$.}

We have the following version of the Morse inequality.

{\bf Theorem 2.5.} (Theorem 6.1 of \cite{Rad2}) {\it Suppose that there exist
only finitely many prime closed geodesics $\{c_j\}_{1\le j\le p}$ on $(M, F)$,
and $0\le a<b\le \infty$ are regular values of the energy functional $E$.
Define for each $q\in\Z$,
\bea
{M}_q(\ol{\Lm}^b,\ol{\Lm}^a)
&=& \sum_{1\le j\le p,\;a<E(c^m_j)<b}\rank{\ol{C}}_q(E, c^m_j ), \nn\\
{b}_q(\ol{\Lm}^{b},\ol{\Lm}^{a})
&=& \rank H_q(\ol{\Lm}^{b},\ol{\Lm}^{a}). \nn\eea
Then there holds }
\bea
M_q(\ol{\Lm}^{b},\ol{\Lm}^{a}) &-& M_{q-1}(\ol{\Lm}^{b},\ol{\Lm}^{a})
    + \cdots +(-1)^{q}M_0(\ol{\Lm}^{b},\ol{\Lm}^{a}) \nn\\
&\ge& b_q(\ol{\Lm}^{b},\ol{\Lm}^{a}) - b_{q-1}(\ol{\Lm}^{b},\ol{\Lm}^{a})
   + \cdots + (-1)^{q}b_0(\ol{\Lm}^{b},\ol{\Lm}^{a}), \lb{2.12}\\
{M}_q(\ol{\Lm}^{b},\ol{\Lm}^{a}) &\ge& {b}_q(\ol{\Lm}^{b},\ol{\Lm}^{a}).\lb{2.13}
\eea

\setcounter{equation}{0}
\section{Proof of main theorems}

In this section, we give the proofs of the main theorems. 
We prove by contradiction,  by \cite{LyF} we suppose
the following holds:

{\bf (C) There is only one prime closed geodesic $c$ on $(M, F)$.}

{\bf Proof of Theorem 1.2.} 
Let $p:\tilde{M}\rightarrow M$ be the universal covering of $M$
and $\tilde{F}=p^\ast(F)$. Then  $(\tilde{M}, \tilde{F})$
is a compact Finsler manifold and it is locally isometric to $(M, F)$. 
In fact, $\tilde M$ is a compact manifold without boundary follows from the 
property for covering spaces (cf. Theorem 6.7 and Corollary 6.8  of 
\cite {GrH} or Section 2.6 of \cite{Spa}): Since $M$ is a compact manifold 
without boundary,  the universal covering  space $\tilde M$ of $M$ exists by
Theorem 6.7 and Corollary 6.8 of \cite{GrH}.
By the definition of covering space (cf. p. 21 of \cite{GrH}) and Theorem 5.8 
of \cite{GrH}, 
 $\forall x\in M$ there exists a open neighborhood $U_x$ of $x$ 
 which is homeomorphic to the open disk $\{p\in\R^n: ||p||<1\}$ in $\R^n$ and $\overline {U_x} $ being compact  such that 
 $p^{-1}(U_x)=\cup_{\alpha\in\pi_1(M)}V_{x, \alpha}$ with 
 $V_{x,\alpha}\cap V_{x, \beta}=\emptyset$ for $\alpha\neq\beta$
 and each $V_{x, \alpha}$ is homeomorphic to $U_x$ for $\alpha\in\pi_1(M)$.
 Since $M$ is compact, we can find finitely many $\{x_i\}_{1\le i\le k}$ such that
 $M=\cup_{1\le i\le k}U_{x_i}$, thus  $\tilde M=\cup_{\alpha\in\pi_1(M)}V_{x_i, \alpha}$.
 Since $\pi_1(M)$ is finite, $\tilde M$ is a compact manifold without boundary.
 
If $\tilde{d}$ is a closed geodesic on $(\tilde{M}, \tilde{F})$,
then $d=p(\tilde{d})$ is a closed geodesic on $(M, F)$ and $P_{\tilde{d}}=P_d$.
In fact,  $P_{\tilde{d}}=P_d$ follows from the fact that 
the Jacobi  field equations along $d$ and $\tilde d$  are the same,
while   $P_\gamma$ for a closed geodesic $\gamma$ can be defined 
by using the Jacobi field along  $\gamma$, cf. Section 3 of \cite{Rad4}.
Since $(M, F)$ is bumpy, $d$ is non-degenerate, i.e. $1\notin\sigma(P_d)$, therefore $\tilde{d}$
is non-degenerate also,  hence $(\tilde{M}, \tilde{F})$ is bumpy.
Now we have the following two cases:

{\bf Case 1.} {\it $H^*(\tilde{M};\Q)\neq T_{d,h+1}(x)$.}

In this case, there must be infinitely many prime closed geodesics
 on $(\tilde{M}, \tilde{F})$ by \cite{ViS} of M. Vigu\'e-Poirrier and D. Sullivan.
 In fact, the Betti numbers of the free loop space $\Lambda \tilde M$  are unbounded and 
  the theorem in \cite{GrM} of Gromoll-Meyer can be applied to the Finsler manifolds as well.
  Hence there must be infinitely many prime closed geodesics on
$(M, F)$ since $\pi_1(M)$ is finite. This proves Theorem 1.2 in this case.

{\bf Case 2.} {\it $H^*(\tilde{M};\Q)= T_{d,h+1}(x)$.}

In this case,  as in Case 1, it is sufficient to consider the case  there are
finitely many prime closed geodesics on $(\tilde M, \tilde F)$.

Denote by $\{\tilde{d_1}, \tilde{d_2},\ldots,\tilde{d_k}\}$
the prime closed geodesics on  $(\tilde{M}, \tilde{F})$,
then we have $k\ge 2$ by \cite{DLW}.
If $\pi_1(M)=0$, then Theorem 1.2  holds since $\tilde M=M$ in this case,  thus it remains to consider the case that $\pi_1(M)\neq 0$. For $1\le i\le k$, we have $p(\tilde{d_i})=c^{m_i}$ for some $m_i\in\N$
by the assumption $(C)$.  By  translating the parameters if necessary,
we may assume $p(\tilde d_j(0))=c^{m_j}(0)=c(0)$ for $1\le j\le k$.

Note that for each $i\in \{2, \ldots, k\}$, there exists a covering transformation
$f_i:(\tilde{M}, \tilde{F})\rightarrow(\tilde{M}, \tilde{F})$
such that $f_i(\tilde{d_i}(0))=\tilde d_1(0)$. By the definition of the Finsler metric $\tilde F$ 
on $\tilde M$, the map $f_i$ is an isometry on $(\tilde M, \tilde F)$. Therefore $f_i(\tilde d_i)$
is a closed geodesic started at $f_i(\tilde{d_i}(0))=\tilde d_1(0)$.
By the property of covering transformation, we have 
\bea p(f_i(\tilde d_i(t)))=p(\tilde d_i(t))=c^{m_i}(t)=c(m_it),\quad \forall t\in\R.\lb{3.1}\eea
On the other hand, 
\bea p(\tilde d_1(t))=c^{m_1}(t)=c(m_1t),\quad \forall t\in\R.\lb{3.2}\eea
 Hence we have $f_i(\tilde d_i(t))=\tilde d_1(\frac{m_i}{m_1}t)$ for $\forall t\in\R$.
 Since $f_i(\tilde d_i)$ is a closed geodesic and $\tilde d_1$ is a prime closed 
 geodesic,  we have $m_1|m_i$.  Exchanging  $\tilde d_1$ and $\tilde d_i$,
 we obtain $m_i|m_1$, and then $m_i=m_1$. This yields
 $f_i(\tilde{d_i})=\tilde{d}_1$.
Since $f_i$ is an isometry on $(\tilde{M}, \tilde{F})$,
it preserves the energy functional, i.e., $E(\gamma)=E(f_i(\gamma))$
for any $\gamma\in\Lambda(\tilde M)$. This implies $i(\tilde{d_i}^m)=i(\tilde{d_1}^m)$
for any $m\in\N$. By Proposition 2.1, we have
\bea
M_q(\tilde M)\equiv M_q(\ol{\Lm\tilde M},\ol{\Lm\tilde M}^0)
&=& \sum_{1\le j\le k,\;m\in\N}\rank{\ol{C}}_q(E, \tilde{d}^m_j ) \nn\\
\nn\\
&=& \sum_{1\le j\le k}\, ^\#\{m\,|\,i(\tilde{d}_j^m)-i(\tilde{d}_j)\in 2\Z\;\;{\it and}\;\;
                   q=i(\tilde{d}_j^m)\}
\nn\\
&=&k\,^\#\{m\,|\,i(\tilde{d}_1^m)-i(\tilde{d}_1)\in 2\Z\;\;{\it and}\;\;
                   q=i(\tilde{d}_1^m)\}
\lb{3.3}\eea
By Bott formula, c.f. \cite{Bot}, $i(\tilde{d_1^m})\ge i(\tilde{d_1})$ for $m\in\N$,
then we have $M_q(\tilde M)=0$ for $q<i(\tilde{d_1})$ by (\ref{3.3}).
By Proposition 2.1 and (\ref{3.3}), we have
$M_{i(\tilde{d_1})}(\tilde M)=k\,^\#\{m\,|\,i(\tilde{d}_1^m)=i(\tilde{d}_1)\}$
and $M_{i(\tilde{d_1})+1}(\tilde M)=0$.
We claim that $i(\tilde{d_1})= d-1$.
In fact, by Proposition 2.3 and Theorem 2.5,
$M_{d-1}(\tilde M)\ge b_{d-1}(\tilde M)=1$, this implies $i(\tilde{d_1})\le d-1$.
By  (\ref{3.3}), there exists $m\in\N$ such that $d-1=i(\tilde d_1^m)$
and $i(\tilde d_1^m)-i(\tilde d_1)\in2\Z$, this implies $i(\tilde d_1)-(d-1)\in2\Z$.
Thus if $i(\tilde{d_1})< d-1$, we must have $i(\tilde{d_1})< d-2$ holds.
Hence by Theorem 2.5 and Proposition 2.3,
\bea &&-k\,^\#\{m\,|\,i(\tilde{d}_1^m)=i(\tilde{d}_1)\}
\nn\\=&&
M_{i(\tilde{d_1})+1} (\tilde M)- M_{i(\tilde{d_1})} (\tilde M)+ \cdots +(-1)^{i(\tilde{d_1})+1}M_0(\tilde M)\nn\\\ge&&
 b_{i(\tilde{d_1})+1} (\tilde M)- b_{i(\tilde{d_1})}(\tilde M)+ \cdots + (-1)^{i(\tilde{d_1})+1}b_0(\tilde M) \nn\\=&&0.\lb{3.4}\eea
This contradiction proves that $i(\tilde{d_1})= d-1$.
By Theorem 2.5 and Proposition 2.3 again,
\bea &&-k\,^\#\{m\,|\,i(\tilde{d}_1^m)=i(\tilde{d}_1)=d-1\}
\nn\\=&&
M_{d}(\tilde M) - M_{d-1} (\tilde M)+ \cdots +(-1)^{d}M_0(\tilde M)\nn\\\ge&&
 b_{d}(\tilde M) - b_{d-1}(\tilde M)+ \cdots + (-1)^{d}b_0(\tilde M) \nn\\=&&-1.\lb{3.5}\eea
This contradicts to $k\ge 2$ and  proves Theorem 1.2.\hfill\hb

{\bf Proof of Theorem 1.3.}  Let $\tilde M$ be the orientable double cover of $M$ if $M$ is 
not orientable and $\tilde M=M$ otherwise.
Let  $p:\tilde{M}\rightarrow M$ be the projection
and $\tilde{F}=p^\ast(F)$. Then  $(\tilde{M}, \tilde{F})$
is a compact orientable Finsler manifold without boundary of dimension 2 as 
proved in Theorem 1.2.
Due to the classification of surfaces, cf. Chapter 9 of \cite{Hir},   $\tilde M$
is either a sphere or a torus of genus $g$ with $g\in\N$.
In fact, by Section 9.1 of \cite{Hir}, if $M$ is a surface with genus $g$,
then its Euler characteristic $\chi(M)=2-2g$ if $M$ is orientable,
while $\chi(M)=2-g$ if $M$ is not orientable; therefore, by Section 5.2 of \cite{Hir}, 
 $\chi(\tilde M)=2-2g$ if $M$ is orientable,
while $\chi(\tilde M)=2\chi(M)=4-2g$ if $M$ is not orientable.

We have the following two cases:

{\bf Case 1.} {\it $\tilde M$ is a torus of genus $g$ with $g\in\N$.
}

In this case, we have $b_1(\tilde M)=2g\ge 2$.  Choose 
two linearly independent generators $\alpha, \beta$ of $H_1(\tilde M)$
and $c_m\in \Lambda\tilde M$  such that $[c_m]=\alpha\beta^m$ and
 $c_m$ being a closed geodesic for $m\in\N$. In fact, since $H_1(\tilde M)=\pi_1(\tilde M)/[\pi_1(\tilde M), \pi_1(\tilde M)]$, one can  choose some $c_m^\prime\in\Lambda \tilde M$
 such that $[c_m^\prime]=\alpha\beta^m$.
 Since the energy functional $E$ satisfies the Palais-Smale condition (cf. Theorem 4.6 of \cite{Mer}),
 the sequence $\{\phi_s(c_m^\prime)\}_{s\ge 0}$ has a limit  point $c_m$, 
 where $\phi_s$ is  the  negative
 gradient flow of $E$.  Thus $c_m$ is a closed geodesic and $c_m$ is homotopic   to $\phi_s(c_m^\prime)$ for some $s$ large enough, and then $[c_m]=[\phi_s(c_m^\prime)]=[c_m^\prime]\in H_1(\tilde M)$.
Clearly $c_m$s are distinct prime  closed geodesics. 
Therefore there are infinitely many prime closed geodesics on $\tilde M$.
Since $p: \tilde M\rightarrow M$ is $k$-fold cover with $k\le 2$, 
there are infinitely many prime closed geodesics on $ M$ also.
This proves Theorem 1.3 in this case.

{\bf Case 2.} {\it $\tilde M=S^2$.
}

As above it is sufficient to consider the case that 
$(S^2, \tilde F)$ possess finitely many prime closed geodesics. 
Denote by $\{\tilde{d_1}, \tilde{d_2},\ldots,\tilde{d_k}\}$
the prime closed geodesics on  $(S^2, \tilde{F})$.
Thus for $1\le i\le k$, we have $p(\tilde{d_i})=c^{m_i}$ for some $m_i\in\N$
by the assumption $(C)$.  As in Theorem 1.2, we have
 $m_i=m_1$ and $i(\tilde{d_i}^m)=i(\tilde{d_1}^m)$
 and $k_l(\tilde{d_i}^m)^{\pm 1}=k_l(\tilde{d_1}^m)^{\pm 1}$
 for any $m\in\N,\,l\in\Z$, where $k_l(\tilde {d_i}^m)^{\pm 1}=\rank H_{q-i(\tilde {d_i}^m)}(N_{\tilde {d_i}^m}^-\cup\{\tilde {d_i}^m\}, N_{\tilde {d_i}^m}^-)^{\pm 1}$. In fact, each covering transformation  $f_i$ is an isometry on $(S^2, \tilde{F})$,
it preserves the energy functional, i.e., $E(\gamma)=E(f_i(\gamma))$
for any $\gamma\in\Lambda(S^2)$, therefore $f_i$ induces a homeomorphism 
between characteristic manifolds at $\tilde{d_1}^m$ and $\tilde{d_i}^m$. Hence $i(\tilde{d_i}^m)=i(\tilde{d_1}^m)$, $\nu(\tilde{d_i}^m)=\nu(\tilde{d_1}^m)$ and $k_l(\tilde{d_i}^m)^{\pm 1}=k_l(\tilde{d_1}^m)^{\pm 1}$ for any $m\in\N$. Therefore we obtain $\hat\chi(\tilde{d_i})=\hat\chi(\tilde d_1)$.
 By Lemma 4.3 and Theorem 4.4 of \cite{LoW},  we have $\hat i(\tilde d_1)>0$.
 By Proposition 2.4 or Theorem 4.4 of \cite{LoW}, we have
\be k\frac{\hat\chi(\tilde d_1)}{\hat{i}(\tilde d_1)}=
\sum_{j=1}^k\frac{\hat\chi(\tilde d_j)}{\hat{i}(\tilde d_j)}=B(2, 1)=1.
\lb{3.4}\ee
Since  $\hat\chi(\tilde d_1)\in\Q$, this implies $\hat{i}(\tilde d_1)\in\Q$ also.
Thus  there is no irrationally elliptic closed geodesics on $( S^2, \tilde F)$
 since such a closed geodesic $d$ must satisfy $\hat i(d)\notin\Q$, cf.  Section 3 of \cite{LoW}.
 This contradicts   to \cite{LoW} which claims there exist  at least two irrationally elliptic
prime closed geodesics on every Finsler 2-sphere possessing only finitely many
prime closed geodesics. This proves Theorem 1.3.\hfill\hb   

\noindent {\bf Acknowledgements.} I would like to sincerely thank Professor H.-B. Rademacher and I. Taimanov for pointing out a mistake in the first version of the paper, especially 
to Professor I. Taimanov for explaining  me the existence  results for the Theorem of Fet.
I would like to sincerely thank my
advisor, Professor Yiming Long, for introducing me to the theory of
closed geodesics and for his valuable help and encouragement
to me in all ways.
\bibliographystyle{abbrv}

\end{document}